\documentclass[12pt]{amsart}
\usepackage{amssymb,amsfonts,amsthm}
\usepackage{srcltx}

\usepackage{amsmath}
\usepackage{amscd}
\usepackage{amsfonts}
\usepackage{color}

\newcommand{\Z}{{\mathbb Z}}

\newcommand{\la}{\langle}
\newcommand{\ra}{\rangle}

\newcommand{\dist}{{\mathrm{dist}}}

\newcommand{\DG}{{\mathcal D}}

\newtheorem{theorem}{Theorem}[section]
\newtheorem{ex}[theorem]{Example}

\newtheorem{con}[theorem]{Conjecture}

\theoremstyle{definition}
\newtheorem{df}[theorem]{Definition}

\newtheorem{rk}[theorem]{Remark}
\newtheorem{prob}[theorem]{Problem}

\newcommand{\topp}{{\bf top}}

\newcommand{\bott}{{\bf bot}}

\newcommand{\iv}{^{-1}}

\begin{document}

\title{On the conjugacy growth functions of groups}
\author{V.~S.~Guba}
\address{V. S. Guba, Department of mathematics\\
    Vologda State University\\ Vologda,  160600\\
    Russia}
\email{victorguba@mail.ru}

\author{M.~V.~Sapir}
\address{M. V. Sapir, Department of mathematics\\
    Vanderbilt University\\ Nashville, TN 37240, USA}
\email{m.sapir@vanderbilt.edu}
\maketitle


\begin{abstract}
To every finitely generated group one can assign the conjugacy
growth function that counts the number of conjugacy classes
intersecting a ball of radius $n$. Results of Ivanov and Osin show
that the conjugacy growth function may be constant even if the
(ordinary) growth function is exponential. The aim of this paper is
to provide conjectures, examples and statements that show that in
``normal" cases, groups with exponential growth functions also have
exponential conjugacy growth functions.
\end{abstract}

\section{Introduction}

The investigation of the growth of conjugacy classes was motivated by counting closed geodesics (up to free homotopy) on complete Riemannian manifolds. In the case of a negative upper bound for  the sectional curvature of a complete Riemannian manifold  $M$, there is only one closed geodesic in each free homotopy class.
Strengthening a result of Sinai \cite{Sinai}, Margulis \cite{Mar, MarS} proved that for compact manifolds of pinched negative curvature and
exponential volume growth ~ $\exp(hn)$, the number of primitive closed geodesics of period $\le t$ is approximately
\begin{equation}\label{0}
\frac{\exp(ht)}{ht}.
\end{equation}
In group theoretic terms this implies that the number of primitive conjugacy classes intersecting the ball of radius $n$ in the Cayley graph of the fundamental group of $M$ (with respect to some finite generating set) is between $\frac{1}{Cn}\exp(hn)$ and $\frac{C}{n}\exp(hn)$ for some constant $C>1$. (Recall that a conjugacy class is {\em primitive} if it does not contain elements that are proper powers.)  Basically it means that two cyclically reduced products of generators are conjugate ``almost" only if they are cyclic shifts of each other. So this result should not be very surprising to those who know the theory of Gromov-hyperbolic and relatively hyperbolic groups. And indeed, this result has been generalized several times to larger classes of manifolds and groups (see, for example,  \cite{Link, CK} and references therein).

\begin{df} Let $G=\la X\ra $ be a group generated by a finite set $X$. For every $n$ let $g_c(n)$ be the number of conjugacy classes of $G$ intersecting the ball of radius $n$ in $G$. The function $g_c(n)$ will be called the conjugacy growth function of $G$.
\end{df}

The definition differs from that in \cite{Mar, MarS} and others because we do not consider only primitive conjugacy classes. In the case of hyperbolic, relatively hyperbolic groups or CAT(0)-groups considered in \cite{Mar, MarS, CK, Link} and other papers these definitions give equivalent functions (in the usual sense of asymptotic group theory) because the number of non-primitive conjugacy classes is small comparing to the number of all conjugacy classes. But, for example, in a torsion group without involutions every element is a square, so there are no primitive conjugacy classes at all.

It is known that the conjugacy growth function for arbitrary groups can differ dramatically from the (ordinary) growth function. S. Ivanov \cite{Olbook} constructed the first example of a finitely generated infinite group with finite number of conjugacy classes, and Osin \cite{Osin}, answering an old group theory question (one of the origins of this question is also from dynamics), constructed a finitely generated infinite group with just two conjugacy classes.  The conjugacy growth functions for these groups are eventually constants while the ordinary growth functions are exponential. Still it is not known how widespread this phenomenon is. For example, there are no examples of finitely presented groups with exponential growth function and subexponential conjugacy growth function. Conjectures, examples and theorems of this paper show that for ``ordinary" groups, exponential growth should imply exponential conjugacy growth.

Our results are intentionally not of the most general form. We just demonstrate the ideas which clearly can be used in more general situations.

\section{Amenable groups}

\begin{con}\label{c1} For every amenable group of exponential growth, the conjugacy growth function is exponential. An amenable group with polynomial conjugacy growth function is virtually nilpotent.
\end{con}

For solvable groups Conjecture \ref{c1} should follow from the proof
of Milnor \cite{Milnor} and Wolf \cite{Wolf} of the fact that
solvable groups of subexponential growth are virtually nilpotent.
They prove that if a solvable group has exponential growth, it
contains a free subsemigroup. One ``only" needs to show that the
free subsemigroups are Frattini-embedded (i.e. two elements of the
subsemigroup are conjugate in the ambient groups if and only if they
are cyclic shifts of each other). One can also use a result of
Kropholler \cite{Kro} characterizing finitely generated solvable
groups which do not have sections which are wreath products of a
cyclic group with $\Z$ and results of Osin \cite{Osin1} about the
uniform growth of solvable groups.

The second part of Conjecture \ref{c1} is of course an analog of the celebrated theorem of Gromov \cite{Gromov} about groups with polynomial growth. In this regard it would be interesting to compute the conjugacy growth functions of Grigorchuk (and similar) groups of subexponential growth \cite{Grig}. It is hard to believe that any of these groups have polynomial conjugacy growth, but these groups served as counterexamples to many other conjectures before.

\subsection{Two examples}

\begin{ex} The Baumslag-Solitar group $BS(1,n)=\la a,b \mid b\iv ab=a^n\ra$\ $(n\ge 2)$ has exponential conjugacy growth function.
\end{ex}

\proof It is easy to see that for numbers $k\neq l$ not divisible by $n$,  the elements
$a^k, a^l$ are not conjugate in $BS(1,n)$. The length of $a^k$ in $BS(1,n)$ is $\Theta(\log k)$\footnote{We use the standard Computer Science notation:  we write $g(x)=\Theta(f(x))$ (resp. $g(x)=O(f(x))$) if $\frac1Cf(x)\le g(x)\le C f(x)$ (resp. $g(x)\le Cf(x)$) for some constant $C>1$ and all sufficiently large $x$.}(indeed, if $k=m_1n^{t_1}+\cdots+m_s$ is the $n$-ary representation of $k$, then $a^k=b^{-t_1}a^{m_1}b^{t_1} b^{-t_2}a^{m_2}b^{t_2}\ldots$, so the length of $a^k$ is at most $2(t_1+t_2+\cdots)+m_1+m_2+\cdots$. Hence the conjugacy growth function is exponential.
\endproof

The next example shows that free subsemigroups are indeed helpful in proving that the conjugacy growth is exponential.

\begin{ex} \label{sinfty} Let $S_\infty$ be the group of all permutations of $\Z$ with finite supports. Then the cyclic group $\Z=\la b\ra $ acts on $S_\infty$
by shift. Let  $G=S_\infty \rtimes {\mathbb Z}$. This group is
clearly generated by $b$ and the involution $a=(1,2)$. The conjugacy
growth of $G$ is exponential.
\end{ex}

\proof Indeed, consider the subsemigroup generated by $b^{2}$ and
$b^{2}a$. Every element in this semigroup has the form
$b^{n}ab^{n_1}a\ldots$ where all $n, n_i$ are positive even integers.
By taking a cyclic shift, we can assume that the word ends with $a$
and $n\ge n_1,n_2,\ldots$. Let us prove that no two of such
elements are conjugate (this will obviously imply that the conjugacy
growth is exponential). Indeed, the element
\begin{equation}\label{1}
b^{n}ab^{n_1}a\ldots b^{n_k}a
\end{equation}
 is equal to
\begin{equation}\label{2}
b^{n+n_1+\cdots+n_k}a(n_1+\cdots+n_k)a(n_2+\cdots+n_k)\ldots a
\end{equation}
 where
$a(m)=(1+m, 2+m)$, so $a=a(0)$. Note that since all $n_i$ are even
and positive, the transpositions $a(n_i+\cdots+n_k)$ are all
independent, so they pairwise commute. Therefore the form (\ref{2})
completely determines the word (\ref{1}). If we conjugate an element
$b^sp$ where $p\in S_\infty$ by $b^m$, the result is $b^sp(m)$ where
$p(m)$ is obtained from $p$ by shifting all numbers in the
permutation by $m$. If we conjugate $b^{s}p$ by $t\in S_\infty$, the
result equals $b^{s}t(s)\iv pt$. Note that in (\ref{2}), the
absolute value of the exponent of $b$ is always bigger than all
numbers involved in the permutation
$a(n_1+\cdots+n_k)a(n_2+\cdots+n_k)\ldots a$. Therefore if $b^{m}t$
($t\in S_\infty, m>0$) is conjugate to (\ref{2}) and not equal to
(\ref{2}), then either some numbers in the support of the
permutation $t$ are negative or some of them are bigger than $m$. In
both cases $b^mt$ is not of the form (\ref{2}). So two different
elements of the form (\ref{2}) are not conjugate in $G$.  It remains
to note that the word length of (\ref{1}) in the alphabet $\{a,b\}$
is at most $n+n_1+\cdots+n_k+k$.
\endproof

Example \ref{sinfty} is also interesting because the number of
conjugacy classes of elements of $S_\infty$ intersecting the ball of
radius $n$ grows subexponentially. In contrast, as we showed above,
for $BS(1,n)$, the growth of conjugacy classes of the kernel of the
natural homomorphism $BS(1,n)\to \Z$ is exponential.

\subsection{Nilpotent groups}

Recall that by a theorem of Bass \cite{Bass} if $G$ is a finitely
generated nilpotent group, and $\gamma_S(n)$ its growth function for
some generating set $S$,  then there are constants $A,B>0$ such that
$An^d\leq\gamma_S(n)\leq Bn^d$ for all $n\geq 1$, where
$d=d(G)=\sum_{h\geq 1}hr_h$, $r_h$ is the torsion-free rank of
$G_h/G_{h+1}$ and $G=G_1\supseteq G_2\supseteq
G_3\supseteq\cdots\supseteq G_h\supseteq\cdots$ is the lower central
series of $G$. In particular, for the Heisenberg group $H=\la
x,y,z\mid [x,y]=z, zx=xz, zy=yz\ra$ with two generators, its growth
function is $\Theta(n^4)$.

\begin{ex} The
conjugacy growth function of $H$ is $\Theta(n^2\log n)$.
\end{ex}

\proof Indeed, every element of $H$ is uniquely expressed in the
form $x^ky^lz^m$. The ball $B$ of radius $4n$ consists of  such
elements with $k\le n, l\le n, m\le n^2$ (up to a big O).
Conjugating this element by $x$ amounts to adding $k$ to $m$ (and
preserving $k,l$), conjugating by $y$ amounts to adding $l$ to $m$.
Therefore if $k\ne 0$ or $l\ne 0$, every element of the form $x^ky^lz^m$ can be
conjugated into an element of the form $x^ky^lz^s$ where $0\le s<
d=\gcd(k,l)$. The number of elements of the form $x^ky^lz^m$ in $B$
where either $k=l=0$, $m\le n^2$ or $kl=0, |m|\le \max\{|k|,|l|\}$ does not
exceed $O(n^2)$, so these elements can be ignored. So we can assume
that $k\ne0$ and $l\ne 0$. Then $k=dk_0$, $l=dl_0$ where $k_0, l_0$
are co-prime numbers with absolute values bounded by $O(n/d)$. The
number of such pairs $(k_0,l_0)$ is at most $Cn^2/d^2$ for some
constant $C$. Therefore for every $d$ between $1$ and $n$, the
number of  triples $(k,l,s)$ with $\gcd(k,l)=d$, $0\le s< d$ is at
most $Cn^2/d$ for some constant $C$. Summing up for all $d$ from 1 to
$n$, we get $Cn^2(1+1/2+\cdots)\le C'n^2\log n$. On the other
hand, for every $n>1$ and $d<n$ consider all elements
$x^{k_0d}y^{l_0d}z^s$ where $0< k_0\le n/d, 0<l_0\le n/d, 0\le s<d$,
$\gcd(k_0, l_0)=1$. Again it is easy to see that conjugating such an
element by any non-trivial element of $H$ produces an element not of
this form. Hence these elements are pairwise non-conjugate. The
length of each of these elements is at most $C_1n$ for some constant
$C_1$. By the well-known number theory result (see, for example,
\cite[Page 354]{HW}), the number of such triples is $\Theta(n^2/d)$.
Hence the total number of conjugacy classes intersecting a ball of
radius $n$ in $H$ is at least $C_2n^2(1+1/2+\cdots+1/n)\ge
C_3n^2\log n$ for some constant $C_3>0$ and large enough $n$.
\endproof

Similarly, one can prove that under the assumptions and notation of
Bass' theorem above the conjugacy growth function $g_c(m)$ does not
exceed $Cm^s$ where $s=\sum_{h\geq 1}r_h$.

\begin{prob} Find more precise estimates for the conjugacy growth functions
of finitely generated nilpotent groups.
\end{prob}

\section{Diagram groups}

Let us recall the definition of a {\em diagram group} (see \cite{GS1,GS2,GS3} for more formal definitions). A (semigroup) {\em diagram} is a planar directed labeled graph tesselated into cells, defined up to an isotopy of the plane. Each diagram $\Delta$ has the top path $\topp(\Delta)$, the bottom path $\bott(\Delta)$, the initial and terminal vertices $\iota(\Delta)$ and $\tau(\Delta)$. These are common vertices of $\topp(\Delta)$ and $\bott(\Delta)$.  The whole diagram is situated between the top and the bottom paths, and every edge of $\Delta$ belongs to a (directed) path in $\Delta$ between $\iota(\Delta)$ and $\tau(\Delta)$. More formally, let $X$ be an alphabet. For every $x\in X$ we define the {\em trivial diagram} $\varepsilon(x)$ which is just an edge labeled by $x$. The top and bottom paths of $\varepsilon(x)$ are equal to $\varepsilon(x)$, $\iota(\varepsilon(x))$ and $\tau(\varepsilon(x))$ are the initial and terminal vertices of the edge. If $u$ and $v$ are words in $X$, a {\em cell} $(u\to v)$ is a planar graph consisting of two directed labeled paths, the top path labeled by $u$ and the bottom path labeled by $v$, connecting the same points $\iota(u\to v)$ and $\tau(u\to v)$. There are three operations that can be applied to diagrams in order to obtain new diagrams.

1. {\bf Addition.} Given two diagrams $\Delta_1$ and $\Delta_2$, one can identify $\tau(\Delta_1)$ with $\iota(\Delta_2)$. The resulting planar graph is again a diagram denoted by $\Delta_1+\Delta_2$, whose top (bottom) path is the concatenation of the top (bottom) paths of $\Delta_1$ and $\Delta_2$. If $u=x_1x_2\ldots x_n$ is a word in $X$, then we denote $\varepsilon(x_1)+\varepsilon(x_2)+\cdots + \varepsilon(x_n)$ (i.e. a simple path labeled by $u$) by $\varepsilon(u)$  and call this diagram also {\em trivial}.

2. {\bf Multiplication.} If the label of the bottom path of $\Delta_2$ coincides with the label of the top path of $\Delta_1$, then we can {\em multiply} $\Delta_1$ and $\Delta_2$, identifying $\bott(\Delta_1)$ with $\topp(\Delta_2)$. The new diagram is denoted by $\Delta_1\circ \Delta_2$. The vertices $\iota(\Delta_1\circ \Delta_2)$ and $\tau(\Delta_1\circ\Delta_2)$ coincide with the corresponding vertices of $\Delta_1, \Delta_2$, $\topp(\Delta_1\circ \Delta_2)=\topp(\Delta_1), \bott(\Delta_1\circ \Delta_2)=\bott(\Delta_2)$.

\begin{center} 
\unitlength=1mm
\special{em:linewidth 0.4pt}
\linethickness{0.4pt}
\begin{picture}(124.41,55.00)
\put(1.00,30.00){\circle*{2.00}}
\put(46.00,30.00){\circle*{2.00}}
\put(1.00,30.00){\line(1,0){45.00}}
\bezier{320}(1.00,30.00)(24.00,55.00)(46.00,30.00)
\bezier{332}(1.00,30.00)(24.00,5.00)(46.00,30.00)
\put(24.00,35.00){\makebox(0,0)[cc]{$\Delta_1$}}
\put(24.00,25.00){\makebox(0,0)[cc]{$\Delta_2$}}
\put(24.00,10.00){\makebox(0,0)[cc]{$\Delta_1\circ\Delta_2$}}
\put(66.00,30.00){\circle*{2.00}}
\put(94.00,30.00){\circle*{2.00}}
\put(123.00,30.00){\circle*{2.00}}
\bezier{164}(66.00,30.00)(80.00,45.00)(94.00,30.00)
\bezier{152}(66.00,30.00)(81.00,17.00)(94.00,30.00)
\bezier{172}(94.00,30.00)(109.00,46.00)(123.00,30.00)
\bezier{168}(94.00,30.00)(110.00,15.00)(123.00,30.00)
\put(80.00,30.00){\makebox(0,0)[cc]{$\Delta_1$}}
\put(109.00,30.00){\makebox(0,0)[cc]{$\Delta_2$}}
\put(94.00,10.00){\makebox(0,0)[cc]{$\Delta_1+\Delta_2$}}
\end{picture}
\end{center}

3. {\bf Inversion.} Given a diagram $\Delta$, we can flip it about a horizontal line obtaining a new diagram $\Delta\iv$ whose top (bottom) path coincides with the bottom (top) path of $\Delta$. Note that a cell can have the form $u\to u$ (i.e. $v$ and $u$ can coincide). In this case, in order to distinguish the top from the bottom, we draw a vertical arrow inside the cell from the top to the bottom. In the inverse cell $(u\to u)\iv$ then the arrow points in the opposite direction. So the inverse cell is not the same as the cell itself.

\begin{df} A diagram over a collection of cells $P$ is any planar graph obtained from the trivial diagrams and cells of $P$ by the operations of addition, multiplication and inversion. If the top path of a diagram $\Delta$ is labeled by a word $u$ and the bottom path is labeled by a word $v$, then we call $\Delta$ a $(u,v)$-diagram over $P$.
\end{df}

Two cells in a diagram form a {\em dipole} if the bottom part of the first cell coincides with the top part of the second cell, and the cells are inverses of each other. In this case, we can obtain a new diagram removing the two cells and replacing them by the top path of the first cell. This operation is called {\em elimination of dipoles}. The new diagram is called {\em equivalent} to the initial one. A diagram is called {\em reduced} if it does not contain dipoles. It is proved in \cite[Theorem 3.17]{GS1} that every diagram is equivalent to a unique reduced diagram.

Now let $P=\{c_1,c_2,\ldots\}$ be a collection of cells. We can view it either as a semigroup presentation as in \cite{GS1} or as the collection of cells of a directed 2-complex obtained by identifying all edges of the cells with the same label as in \cite{GS2}. We shall use the first point of view here. If the top and the bottom paths of a diagram are labeled by the same word $u$, we call it a {\em spherical} $(u,u)$-diagram. The diagram group $\DG(P,u)$ corresponding to the collection of cells $P$ and a word $u$ consists of all reduced spherical $(u,u)$-diagrams obtained from these cells and trivial diagrams by using the three operations mentioned above. The product $\Delta_1\Delta_2$ of two  diagrams $\Delta_1$ and $\Delta_2$ is the reduced diagram obtained by removing all dipoles from $\Delta_1\circ\Delta_2$. The fact that $\DG(P,u)$ is a group is proved in \cite{GS1}.

Examples: (Both examples can be found in \cite{GS1}.) If $X$
consists of one letter $x$ and $P$ consists of one cell $x\to x^2$,
then the group $\DG(P,x)$ is the R. Thompson group $F$. If $X$
consists of three letters $a,b,c$ and $P$ consists of three cells
$ab\to a, b\to b, bc\to c$, then the diagram group $\DG(P,ac)$ is
isomorphic to the wreath product $\Z \wr \Z$ \cite{GS3}.

Here are the diagrams representing the two standard generators $x_0,
x_1$ of the R. Thompson group $F$. All edges are labeled by $x$ and
oriented from left to right, so we omit the labels and orientation
of edges.

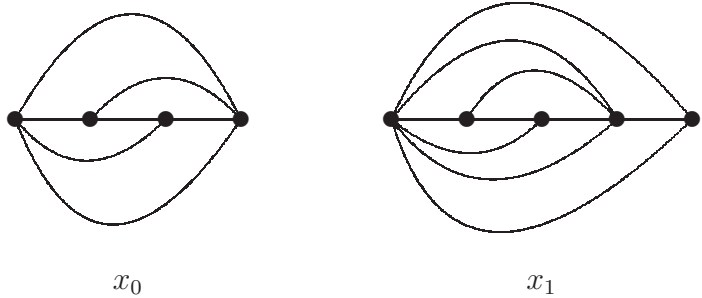
\begin{center} 
\unitlength=1mm
\special{em:linewidth 0.4pt}
\linethickness{0.4pt}
\begin{picture}(94.00,50.00)
\put(3.00,24.00){\circle*{2.00}}
\put(13.00,24.00){\circle*{2.00}}
\put(23.00,24.00){\circle*{2.00}}
\put(33.00,24.00){\circle*{2.00}}
\put(53.00,24.00){\circle*{2.00}}
\put(63.00,24.00){\circle*{2.00}}
\put(73.00,24.00){\circle*{2.00}}
\put(83.00,24.00){\circle*{2.00}}
\put(93.00,24.00){\circle*{2.00}}
\put(3.00,24.00){\line(1,0){10.00}}
\put(13.00,24.00){\line(1,0){10.00}}
\put(23.00,24.00){\line(1,0){10.00}}
\put(53.00,24.00){\line(1,0){10.00}}
\put(63.00,24.00){\line(1,0){10.00}}
\put(73.00,24.00){\line(1,0){10.00}}
\put(83.00,24.00){\line(1,0){10.00}}
\bezier{120}(13.00,24.00)(23.00,35.00)(33.00,24.00)
\bezier{120}(3.00,24.00)(12.00,13.00)(23.00,24.00)
\bezier{256}(3.00,24.00)(19.00,52.00)(33.00,24.00)
\bezier{256}(3.00,24.00)(14.00,-4.00)(33.00,24.00)
\bezier{132}(63.00,24.00)(71.00,37.00)(83.00,24.00)
\bezier{108}(53.00,24.00)(64.00,15.00)(73.00,24.00)
\bezier{208}(53.00,24.00)(70.00,45.00)(83.00,24.00)
\bezier{176}(53.00,24.00)(65.00,8.00)(83.00,24.00)
\bezier{296}(53.00,24.00)(67.00,55.00)(93.00,24.00)
\bezier{296}(53.00,24.00)(62.00,-6.00)(93.00,24.00)
\put(18.00,2.00){\makebox(0,0)[cc]{$x_0$}}
\put(73.00,2.00){\makebox(0,0)[cc]{$x_1$}}
\end{picture}
\end{center}

It is easy to represent, say, $x_0$ as a product of sums of cells and trivial diagrams:
$$x_0=(x\to x^2)\circ (\varepsilon(x)+(x\to x^2))\circ ((x\to x^2)\iv +\varepsilon(x))\circ ((x\to x^2)\iv).$$

There is a natural {\em diagram metric} on every diagram group
$\DG(P,u)$: $\dist(\Delta,\Delta')$ is the number of cells in the
diagram $\Delta\iv\Delta'$. For some finitely generated diagram
groups (such as $F$ or $\Z \wr\Z$) this metric is quasi-isometric to
the word metric \cite{AGS}. In this case, we say that the group
satisfies property B (after J. Burillo). We do not know whether
every finitely generated diagram group satisfies B (see
\cite[Question 1.5]{AGS}).

We shall need the following description of the conjugacy relation in diagram groups.

Let $Q$ be a collection of cells. A spherical $(u,u)$-diagram $\Delta$ over $Q$ is called {\em absolutely reduced} if all its $\circ$-powers $\Delta\circ \Delta\circ \cdots\circ \Delta$ are reduced. Every absolutely reduced diagram $\Delta$ is canonically decomposed as a sum $\Delta_1+\Delta_2+\cdots +\Delta_n$ where each $\Delta_i$ is either a trivial diagram $\varepsilon(u_i)$ or a spherical $(u_i,u_i)$-diagram which is further indecomposable as a sum of spherical diagrams. Indecomposable absolutely reduced spherical diagrams are called {\em simple}.

\begin{theorem}\label{conj} (see \cite[Lemmas 15.14, 15.15, 15.20]{GS1}) (i) Every spherical $(u,u)$-diagram is conjugate to an absolutely reduced spherical $(v,v)$-diagram.

(ii) Suppose that two absolutely reduced diagrams $A$ and $B$ have canonical decompositions $A_1+\cdots+A_m$ and $B_1+\cdots+B_n$ (where $A_i$ is a $(u_i,u_i)$-diagram, $B_j$ is a $(v_j,v_j)$-diagram). Suppose further that $A$ and $B$ are conjugate. Then $m=n$, and $A_i$ is conjugate to $B_i$, that is $A_i=\Gamma_i\iv B_i\Gamma_i$ for some $(v_i,u_i)$-diagram $\Gamma_i$,  $i=1,\ldots,m$.

(iii) If two simple diagrams $A$, $B$ are conjugate then they have the same number of cells. Two trivial diagrams $\varepsilon(u)$ and $\varepsilon(v)$ are conjugate, if and only if $u=v$ modulo the semigroup presentation $P$.
\end{theorem}

\begin{con} The following conditions for a finitely generated diagram group $G$ are equivalent.
\begin{enumerate}
\item $G$ contains a non-Abelian free subsemigroup.
\item The growth function of $G$ is exponential.
\item The conjugacy growth function of $G$ is exponential.
\end{enumerate}
\end{con}

It is clear that $(3)\to (2)$, $(1)\to (2)$. The implication $(2)\to (1)$ (and hence equivalence of these two conditions) is proved as follows. Let $G$ be a finitely
generated diagram group of exponential growth. By \cite[Theorem 5.7]{GS5}, $G$ is embedded into a (split)
extension of a right angled Artin group $A$ by the R. Thompson group $F$. A subgroup of a right angled Artin group is either Abelian, or contains a non-Abelian free subgroup \cite[Lemma 7.6]{GS5}. Hence we can assume that $G\cap A$ is Abelian. A subgroup of $F$ either contains a copy of $\Z\wr  \Z$ or is Abelian \cite[Theorem 21]{GS4}. Since $\Z\wr\Z$ contains a non-Abelian free semigroup, we can assume that the the subgroup $G/(G\cap A)$ of $F$ is Abelian. Hence $G$ is solvable of class at most 2. Since $G$ has exponential growth, it is not virtually nilpotent. By the result of Milnor and Wolf cited above $G$ contains a non-Abelian free subsemigroup.

The next theorem proves a weaker statement than the implication $(1)\to (3)$.

\begin{theorem} Every finitely generated diagram group containing the wreath product $\Z \wr \Z$ (in particular, the R.Thompson group $F$ \cite{GS2}) has exponential conjugacy growth function.
\end{theorem}

\proof Consider first the group $\Z\wr\Z$ itself, that is the diagram group $\DG(P,ac)$ for  $P=\{ab\to a$, $b\to b$, $bc\to c\}$. Let $\pi$ be the cell $b\to b$, and $n_0,\ldots, n_k$ be positive integers. Let $\Delta(n_0,\ldots,n_k)$ be the following diagram:
\begin{equation}\label{5}
\varepsilon(a)+\pi^{n_0}+\cdots+\pi^{n_k}+\varepsilon(c).
\end{equation}
Let $n=n_0+\cdots+n_k$, $\Gamma(n)$ be the diagram obtained as the
following product $((ab,a)\iv+\varepsilon(c))\circ((ab,a)\iv
+\varepsilon(bc))\circ \cdots \circ ((ab,a)\iv+\varepsilon(b^kc)$.
It is clearly an $(ac,ab^kc)$-diagram. Finally let
$$A(n_0,\ldots,n_k)=\Gamma(n) \Delta(n_0,\ldots,n_k)\Gamma(n)\iv.$$
This is a spherical $(ac,ac)$-diagram from $\DG(P,ac)$.

Note that each component $\pi^{n_i}$ is absolutely reduced, hence
simple. Thus the decomposition (\ref{5}) is the canonical
decomposition of the absolutely reduced diagram
$\Delta(n_0,\ldots,n_k)$. Since $A(n_0,\ldots,n_k)$ is a conjugate
of $\Delta(n_0,\ldots,n_k)$, two diagrams $A(n_0,\ldots,n_k)$ and
$A(n_0',\ldots,n_l')$ are conjugate in $\DG(Q,u)$ if and only if the
diagrams $\Delta(n_0,\ldots,n_k)$ and $\Delta'(n_0',\ldots,n'_l)$
are conjugate. By Theorem \ref{conj} then $k=l$ and
$n_0=n_0',\ldots, n_k=n_k'$. Note also that each diagram
$A(n_0,\ldots,n_k)$ has $n+2(k+1)\le 3n$ cells, whence by property
B, it belongs to the ball of radius $O(n)$ in $\Z\wr \Z$. The number
of solutions $n=n_0+\cdots+n_k$ in positive integers is equal to
$2^{n-1}$ ($k\ge 0$). Hence the conjugacy growth of $\Z\wr\Z$ is
exponential.

Now suppose that for some collection of cells (semigroup
presentation) $Q$ and some word $u$ we have $\DG(Q,u)\ge\Z\wr\Z$. We
shall use the fact that $\Z\wr\Z$ is a {\em rigid} diagram group in
terminology of \cite{GS4,GS3}: by \cite[Theorem 24]{GS4} (for
another formulation see \cite[Section 10]{GS3}), there exists an
embedding $\Psi$ of $\Z\wr \Z$ into $\DG(Q,u)$ {\em induced} by a
$(\psi(ac),u)$-diagram $\Gamma$, and a map $\psi$ that takes letters
$a,b,c$ to words $\psi(a), \psi(b), \psi(c)$ over the alphabet of
$Q$, and each of the three cells $x\to y$ of $P$ to a non-trivial
$(\psi(x),\psi(y))$-diagram $\psi(x\to y)$ over $Q$. The map $\Psi$
takes each $(ac,ac)$-diagram $\Delta$ of $\DG(P,ac)$ to the diagram
$\Gamma\iv \psi(\Delta)\Gamma$ where $\psi(\Delta)$ is obtained from
$\Delta$ by replacing every edge $\varepsilon(e)$ by the path
$\varepsilon(\psi(e))$ and every cell $\pi$ by the diagram
$\psi(\pi)$.

Note that the number of cells in $\Psi(\Delta)$ is big-O of the number of cells in $\Delta$. Hence if, as above,  we let $n=n_0+\cdots+n_k$, then by property B, the diagram $\Psi(A(n_0,\ldots,n_k))$ is in the ball of radius $O(n)$. Also note that $\Psi(A(n_0,\ldots,n_k))$ is conjugate to the diagram $\varepsilon(\psi(a))+\psi(\pi)^{n_0}+\cdots+\psi(\pi)^{n_0}+\varepsilon(\psi(c))$. By Theorem \ref{conj}, Part (i), $\psi(\pi)$ is conjugate to an absolutely reduced diagram that is a sum of simple or trivial diagrams $\Pi_1+\cdots+\Pi_s$, and the number of simple diagrams among $\Pi_i$ is not zero. Therefore $\Psi(A(n_0,\ldots, n_k)$ is conjugate to the diagram $$\varepsilon(\psi(a))+(\Pi_1^{n_0}+\cdots+\Pi_s^{n_0})+\cdots+(\Pi_1^{n_k}+\cdots+\Pi_s^{n_k})+\varepsilon(\psi(c)).$$
By Theorem \ref{conj}, Parts (ii) and  (iii), different diagrams of the form $\Psi(A(n_0,\ldots,n_k))$ are not conjugate  in $\DG(Q,u)$. This proves that the conjugacy growth of $\DG(Q,u)$ is exponential.
\endproof

\section{Groups acting on simplicial trees}
\begin{con} Suppose that $G$ acts on a simplicial tree non-trivially and faithfully. Then the conjugacy growth function of $G$ is exponential provided the growth function of $G$ is exponential.
\end{con}

\begin{theorem}\label{th1} Let $G$ be the HNN extension of a group $H$ with associated subgroups $A,B$ such that $AB \cup BA\ne H$. Then the conjugacy growth function of $G$ is exponential.
\end{theorem}

\begin{rk} Note that the Baumslag-Solitar group $BS(m,n)=\la a,b\mid b\iv
a^mb=a^n\ra$, satisfies the conditions of Theorem \ref{th1} if and
only if $m,n$ are not co-prime. There is no doubt that the
conclusion of the theorem remains true for every $m,n$ except
$m=n=1$ (in which case $BS(m,n)=\Z^2$ and the conjugacy growth
function is quadratic).
\end{rk}

\proof Let $a\in G\setminus (AB\cup BA)$, $t$ be the free letter of the $HNN$-extension. Consider the subsemigroup $S$ generated by $t, t a$. Every word in $S$ has the form
\begin{equation}\label{3}
t^{n_1}at^{n_2}\ldots at^{n_{k+1}}
\end{equation}
were all $n_i\ge 0$, and $n_1,\ldots,n_{k}>0$.

Taking a cyclic shift if necessary, we will assume that $n_{k+1}=0$.

We fix the presentation of $G$ which consists of all relations of
$H$ plus the conjugacy relations $ut=tv$ of the HNN-extension (here
$u\in A, v\in B$). Suppose that two words of the form (\ref{3}) are
conjugate in $G$. Then there exists an annular (Schupp) diagram
$\Delta$ for this conjugacy with the inner (outer) label of the form
$t^{n_1}at^{n_2}\ldots t^{n_k}a$ (resp.  $t^{m_1}at^{m_2}\ldots
t^{m_l}a$) (see \cite{MS}). We can assume that $\Delta$ has minimal
possible number of cells. Since $G$ is an HNN-extension, following
Miller and Schupp \cite{MS}, we can consider the $t$-bands, that is
sequences of cells corresponding to the HNN-relations such that
every two consecutive cells have a $t$-edge in common. Since our
HNN-relations include all relations of the form $ut=tv$, $u\in A,
v\in B$, and since $\Delta$ has minimal possible number of cells,
every $t$-band must have length $1$ or $0$. It is proved in
\cite{MS}, that $\Delta$ cannot contain $t$-annuli (i.e. $t$-bands
with same first and last $t$-edges). Also a $t$-band cannot connect
two $t$-edges on the same boundary component of $\Delta$ because all
$n_i\ge 0$. Hence every $t$-band connects a $t$-edge of the inner
boundary of $\Delta$ with a $t$-edge on the outer boundary of
$\Delta$. Therefore the $t$-bands subdivide $\Delta$ into a number
of disc van Kampen diagrams $\Gamma_1,\ldots,\Gamma_s$ without
HNN-cells, hence without $t$-edges. Thus each $\Gamma_i$ is a
diagram over $H$, and so it consists of one cell corresponding to a
relation of $H$. Note that the boundary of $\Gamma_i$ is of the form
$p_1(\Gamma_i)\iv q_1(\Gamma_i)p_2(\Gamma_i) q_2(\Gamma_i)$ (read
counterclockwise) where $p_1(\Gamma_i)$ is a subpath of the inner
boundary of $\Delta$, $p_2(\Gamma_i)$ is a subpath of the outer
boundary of $\Delta$, $q_1(\Gamma_i)$ and $q_2(\Gamma_i)$ are sides
of two consecutive $t$-bands.

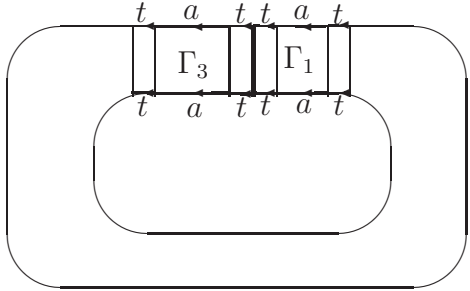
\begin{center}
\unitlength .6mm 
\linethickness{0.4pt}
\ifx\plotpoint\undefined\newsavebox{\plotpoint}\fi 
\begin{picture}(112.25,66.75)(0,0)
\put(61.375,34.75){\oval(101.75,58)[]}
\put(62.625,33.25){\oval(65.75,31)[]}
\put(38.25,63.5){\line(0,-1){15.25}}
\put(81.5,63.5){\line(0,-1){15.25}}
\put(59.75,63.5){\line(0,-1){15.25}}
\put(65.25,63.5){\line(0,-1){15.25}}
\put(43.25,63.5){\line(0,-1){15.25}}
\put(86.5,63.5){\line(0,-1){15.25}}
\put(64.75,63.5){\line(0,-1){15.25}}
\put(70.25,63.5){\line(0,-1){15.25}}
\put(40.75,63.75){\vector(-1,0){.07}}\put(43,63.75){\line(-1,0){4.5}}
\put(40.625,49){\vector(-1,0){.07}}\put(43.25,49){\line(-1,0){5.25}}
\put(51.5,63.5){\vector(-1,0){.07}}\put(59.75,63.5){\line(-1,0){16.5}}
\put(51.625,48.875){\vector(-1,0){.07}}\multiput(59.75,49)(-2.03125,-.03125){8}{\line(-1,0){2.03125}}
\put(62.375,63.75){\vector(-1,0){.07}}\put(65,63.75){\line(-1,0){5.25}}
\put(62.125,48.625){\vector(-1,0){.07}}\multiput(64.75,48.5)(-.65625,.03125){8}{\line(-1,0){.65625}}
\put(67.5,63.75){\vector(-1,0){.07}}\multiput(70,63.5)(-.3333333,.0333333){15}{\line(-1,0){.3333333}}
\put(67.5,48.875){\vector(-1,0){.07}}\multiput(69.75,48.75)(-.5625,.03125){8}{\line(-1,0){.5625}}
\put(75.625,63.625){\vector(-1,0){.07}}\multiput(81,63.5)(-1.34375,.03125){8}{\line(-1,0){1.34375}}
\put(75.75,48.875){\vector(-1,0){.07}}\multiput(81.25,49)(-1.375,-.03125){8}{\line(-1,0){1.375}}
\put(83.875,63.875){\vector(-1,0){.07}}\multiput(86,63.75)(-.53125,.03125){8}{\line(-1,0){.53125}}
\put(84,48.875){\vector(-1,0){.07}}\multiput(86.5,49)(-.625,-.03125){8}{\line(-1,0){.625}}
\put(40.25,66.75){\makebox(0,0)[cc]{$t$}}
\put(51.25,66.75){\makebox(0,0)[cc]{$a$}}
\put(62.25,66.5){\makebox(0,0)[cc]{$t$}}
\put(67.75,66.5){\makebox(0,0)[cc]{$t$}}
\put(75.75,66.75){\makebox(0,0)[cc]{$a$}}
\put(83.75,66.25){\makebox(0,0)[cc]{$t$}}
\put(40.5,46){\makebox(0,0)[cc]{$t$}}
\put(52,45){\makebox(0,0)[cc]{$a$}}
\put(62.5,45.75){\makebox(0,0)[cc]{$t$}}
\put(67.75,46.25){\makebox(0,0)[cc]{$t$}}
\put(76.25,46){\makebox(0,0)[cc]{$a$}}
\put(84,46.25){\makebox(0,0)[cc]{$t$}}
\put(51.75,55.25){\makebox(0,0)[cc]{$\Gamma_3$}}
\put(75.25,56.25){\makebox(0,0)[cc]{$\Gamma_1$}}
\end{picture}

\end{center}

Since $p_j(\Gamma_i)$ does not contain $t$-edges, we obtain that either $p_j(\Gamma_i)$ is empty or it is an $a$-edge, $j=1,2$. Since all $n_1,\ldots,n_k$ are non-negative, we have that the label $u_i$ of the path $q_1(\Gamma_i)$ is from $A$, the label $v_i$ of the path $q_2(\Gamma_i)$ is from $B$. Hence each diagram $\Gamma_i$ corresponds to a relation of $H$ of the form $a^{-\varepsilon_i} u_ia^{\delta_i}v_i$ where $\varepsilon_i, \delta_i\in \{0,1\}$. Since $a\not\in AB\cup BA$, we conclude that for every $i$, $\varepsilon_i=\delta_i$. Therefore the $t$-bands of $\Delta$ provide a correspondence between letters of the words written on the inner and outer boundaries of $\Delta$. This correspondence shows that these two words are cyclic shifts of each other, so $l=k$ and for some $j$ we have $n_1=m_j, n_2=m_{j+1}, \ldots$ (addition in the indices is modulo $k$). Thus we prove that two words of the form (\ref{3}) are conjugate in $G$ if and only if they are cyclic shifts of each other. Hence the conjugacy growth function of $G$ is exponential.
\endproof

\end{document}